%% file: manuscript.tex
\documentclass[%
11pt,
onecolumn,
tightenlines,
superscriptaddress,
preprintnumbers,
nofootinbib,
amsmath,amssymb,amsthm,
physrev,
eqsecnum,
]{revtex4-1}

\input{preamble}

\newcommand{\bfFr}{\bfF_{relax}}

\begin{document}

\preprint{To appear in Journal of the Mechanics and Physics of Solids (\url{https://doi.org/10.1016/j.jmps.2021.104499})}

\title{\Large{Surface Growth in Deformable Solids using an Eulerian Formulation}}

\author{Kiana Naghibzadeh}
    \email{snaghibz@andrew.cmu.edu}
    \affiliation{Department of Civil and Environmental Engineering, Carnegie Mellon University}

\author{Noel Walkington}
    \affiliation{Center for Nonlinear Analysis, Department of Mathematical Sciences, Carnegie Mellon University}
    
\author{Kaushik Dayal}
    \affiliation{Department of Civil and Environmental Engineering, Carnegie Mellon University}
    \affiliation{Center for Nonlinear Analysis, Department of Mathematical Sciences, Carnegie Mellon University}
    \affiliation{Department of Materials Science and Engineering, Carnegie Mellon University}
    
\date{\today}

\begin{abstract}
    Growth occurs in a wide range of systems ranging from biological tissue to additive manufacturing.
    This work considers surface growth, in which mass is added to the boundary of a continuum body from the ambient medium or from within the body.
    In contrast to bulk growth in the interior, the description of surface growth requires the addition of new continuum particles to the body.
    This is challenging for standard continuum formulations for solids that are meant for situations with a fixed amount of material.
    Recent approaches to handle this have used, for instance, higher-dimensional time-evolving reference configurations.
    
    In this work, an Eulerian approach to this problem is formulated, enabling the side-stepping of the issue of constructing the reference configuration.
    However, this raises the complementary challenge of determining the stress response of the solid, which typically requires the deformation gradient that is not immediately available in the Eulerian formulation.
    To resolve this, the approach introduces additional kinematic descriptors, namely the relaxed zero-stress deformation and the elastic deformation; in contrast to the deformation gradient, these have the important advantage that they are not required to satisfy kinematic compatibility.
    The zero-stress deformation and the elastic deformation are used to eliminate the deformation gradient from the formulation, with the evolution of the elastic deformation shown to be governed by a transport equation.
    The resulting model has only the density, velocity, and elastic deformation as variables in the Eulerian setting.
    The proposed method is applied to simplified examples that demonstrate non-normal growth and growth with boundary tractions.

    The introduction in this formulation of the relaxed deformation and the elastic deformation provides a description of surface growth whereby the added material can bring in its own kinematic information.
    Loosely, the added material ``brings in its own reference configuration'' through the specification of the relaxed deformation and the elastic deformation of the added material.
    This kinematic description enables, e.g., modeling of non-normal growth using a standard normal growth velocity and a simple approach to prescribing boundary conditions.

\end{abstract}

\maketitle


\section{Introduction} \label{sec:intro}

Growth, i.e., the addition of material, occurs in a wide variety of systems ranging from biological tissue to additive manufacturing. 
It can be divided into bulk growth and surface growth\footnote{
    In the biological context, these are sometimes called interstitial and appositional growth, respectively \cite{ateshian2007theory}.
}.
Bulk growth is typically studied by assuming that the set of the material points is unchanged, and that growth occurs through a change in the \textit{referential} density due to a distributed mass source.
Because the set of material particles is constant, the standard framework of continuum mechanics is largely applicable with minimal changes; examples of different approaches include the use of growth tensors to define the anelastic shape change \cite{rodriguez1994stress, chenchiah2014energy, reina2017incompressible}, geometrical approaches \cite{yavari2010geometric}, and mixture theory \cite{humphrey2002constrained}.

In contrast, surface growth involves the accretion and ablation of mass at the boundary.
This requires the introduction of new material particles to the body \cite{skalak1982analytical}, and poses interesting challenges.
Surface growth occurs in a range of settings including in biological tissue \cite{taber1995biomechanics}, planetary formation \cite{brown1963gravitational, kadish2008stresses}, solidification and casting processes \cite{schwerdtfeger1998stress}, etching process in silicon wafers \cite{rao2000modelling}, additive manufacturing processes \cite{drozdov1998continuous}, construction of masonry structures \cite{bacigalupo2012effects}, and cell motility via assembly and disassembly of Actin networks \cite{papadopoulos2010surface}.

The addition and removal of material points in surface growth causes difficulties in the usual procedure of defining a fixed set of points as the reference configuration. 
Dealing with this in linear elasticity is somewhat easier \cite{ong2004equations, bacigalupo2012effects,naumov1994mechanics,kadish2005stresses}; for instance, \cite{brown1963gravitational} is an early study of surface growth in linear elasticity, and they showed that the construction of a body by accretion can cause residual stresses which depends on the history of the construction of the body.

However, in nonlinear elastic descriptions of surface growth, the definition of the reference configuration is more subtle.
Pioneering work by Skalak, Hoger, and co-workers provided a systematic approach to study the kinematics of surface growth \cite{skalak1982analytical,skalak1997kinematics}. 
They introduced the \textit{time of attachment} $\tau$ as an additional descriptor to account for the history of attachment. 
However, these studies did not consider deformation, and assumed that the body was rigid. 
Recent work has extended these ideas to account for deformation and stress \cite{sozio2017nonlinear, sozio2020nonlinear,sozio2019nonlinear,tomassetti2016steady,abi2019kinetics,abi2020surface,abeyaratne2020treadmilling}.
Specifically, Cohen, Abeyaratne, and co-workers \cite{,tomassetti2016steady,abi2019kinetics,abi2020surface,abeyaratne2020treadmilling} use a 4-dimensional hypersurface manifold as the reference configuration, in which the 4th dimension accounts for the attachment time $\tau$ of each particle. 
These studies all use a Lagrangian formulation, i.e., roughly the kinematic description is based on reference particles in an evolving reference configuration.
Building on this body of work, they recently proposed a numerical method for the Lagrangian formulation \cite{von2020morphogenesis}.
They solve the equations in an evolving reference configuration numerically, with its shape computed using a volume-conserved regularization, and the new reference configuration requiring re-meshing as it evolves.

In this work, we develop an Eulerian description of surface growth.
This has the advantage that we do not need to explicitly calculate the reference configuration, e.g. \cite{clayton2013nonlinear,clayton2019nonlinear,clayton2020aplace}.
On the other hand, in a standard solid, the stress response requires knowledge of the deformation gradient $\bfF$, and the latter is not straightforward to calculate in an Eulerian setting.
Prior work, notably \cite{liu2001eulerian,kamrin2009eulerian,kamrin2012reference}, have developed approaches to find the deformation gradient in an Eulerian setting; specifically, \cite{liu2001eulerian} uses an evolution equation for $\bfF$ that does not require explicit calculation of the reference configuration.

However, while these approaches work well for the standard setting of a fixed set of material particles without growth, there is the challenge of defining the reference state and kinematic information of the added material particles during surface growth.
We therefore introduce two additional kinematic descriptors: the relaxed zero-stress deformation $\bfF_{relax}$ and the elastic deformation $\bfF_e$.
These are related to the deformation gradient through the relation $\bfF_e = \bfF \bfF_{relax}^{-1}$.
We eliminate $\bfF$ -- as well as $\bfF_{relax}$, using that it is transported with the material particles -- and thereby formulate a model that is posed in terms of $\bfF_e$ as the sole kinematic descriptor for the stress response.
Having $\bfF_e$ as the sole kinematic variable provides some important advantages:
(1) the boundary conditions for the added material is specified in terms of $\bfF_e$ for the added particles, which can be obtained in a transparent manner from the given stress state of the added particles;
(2) $\bfF_e$ has no restrictions of kinematic compatibility, unlike $\bfF$; this is most important in the specifications of boundary conditions of added material, enabling us to specify $\bfF_e$ \textit{solely} from the stress state of the added particles without any concern that it might conflict with the restrictions of kinematic compatibility;
(3) the ability to model non-normal growth using an appropriate $\bfF_e$ in combination with a standard normal growth velocity, potentially providing advantages for numerical methods for free boundary problems that are typically posed with only a normal velocity.
Roughly, we can think of $\bfF_e$ and $\bfFr$ as providing a mechanism for the added material to bring in its own kinematic information to define its own reference configuration.

The evolution of $\bfF_e$ is governed by a transport equation.
The structure of the equation allows us to transparently handle both accretion and ablation: accretion corresponds to inflow, and requires the specification of boundary conditions, while ablation corresponds to outflow and requires no boundary conditions.

Finally, we highlight that the overall structure of surface growth has motivated an Eulerian approach in various prior works.
For instance, \cite{ateshian2007theory} modeled both surface and bulk growth using mixture theory, with the balance equations written in the Eulerian form; however, to compute $\bfF$, he defined the motion of the new particles as a summation of motion of the nearest material point on the boundary and the normal growth velocity, and other model-specific assumptions.
Another body of important contributions in the Eulerian setting is \cite{ganghoffer2011mechanics, ganghoffer2010mechanical, ganghoffer2018combined}; 
similar to these studies, we will use the spatial description to write the balance laws, but compute the deformation using the time rate of change relative to a reference configuration that is not explicitly required to be calculated.

\paragraph*{Organization.}
Section \ref{sec:balance} summarizes the balance laws and jump / boundary conditions, using the view of surface growth as a localized source on the growing surface.
Section \ref{sec:kinematics} discusses the kinematics of the deformation in the Eulerian setting, focusing on introducing the zero-stress deformation and the elastic deformation, the evolution equations for these quantities, and other kinematic issues.
Section \ref{sec:summary}  summarizes the field equations corresponding to our approach.
Section \ref{sec:non-normal-growth} provides an example of non-normal growth; Section \ref{sec:compatible-deformations} provides an example of compatible deformations in a setting with varying elastic deformation; and Section \ref{sec:shear} provides an example of growth accounting for momentum transfer and shear at the boundary.

\paragraph*{Notation.}

The motion is denoted by $\bfx = \bfchi(\bfX , t)$, where $\bfX$ is the location in the reference configuration of a point that is currently at the spatial location $\bfx$.
Also, the inverse of the motion $\bfchi^{-1} (\bfx , t)$  maps the current spatial location of a particle at $\bfx$ to the location $\bfX$ in the reference configuration. 
The deformation gradient is $\bfF = \parderiv{\bfchi (\bfX , t)}{\bfX}$.

Because most of our work is in the Eulerian setting, all differential operators (e.g. $\nabla, \divergence$) imply derivatives with respect to $\bfx$, except where explicitly stated.
Similarly, the argument of various field quantities will implicitly be $\bfx$, except where explicitly stated.

\section{Defining Surface Growth: Surface Sources, Balance Laws, and Constitutive Response}
\label{sec:balance}

We model the process of surface growth through the introduction of sources that are localized on the growing surface.
The source terms can be understood as a coarse-grained approach to treat the complex process of growth without considering the fine details of the processes in the ambient environment outside the growing body; the growth is defined only in terms of net mass and linear momentum transfer\footnote
{
    In addition, we must specify the stress response and initial stress state of the added material; we discuss this in Section \ref{sec:kinematics}.
}.

To define the source terms, we specify the following properties of the added material: the rate of mass addition per unit area, $M$; and the rate of linear momentum addition per unit area, $P$.
From the quantities, we can infer the velocity (or specific momentum) of the added material at the instant of attachment, $\bfv_a := P / M$; and the mass density, $\rho := \frac{M}{\bfv_a \cdot \hat\bfn}$.
Accretion is modeled by $M>0$ and ablation by $M<0$.

We assume that the added particles do not carry angular momentum such as due to individual particle spins.
Therefore, the balance of angular momentum provides only the standard result that the Cauchy stress must be symmetric.

\subsection{Bulk Balance Laws}

The sources corresponding to surface growth are localized on the surface of the body and appear only in the jump conditions.
Therefore, the bulk balance laws are standard:
\begin{align}
    \label{eqn:bulk-balance-mass}
    \text{Mass: } & \parderiv{\rho}{t} + \divergence(\rho \bfv) = 0 
    \\
    \label{eqn:bulk-balance-momentum}
    \text{Momentum: } & \parderiv{\ }{t} (\rho \bfv) + \divergence (\rho \bfv \otimes \bfv)= \rho \bfb + \divergence(\bfsigma)
\end{align}
where $\bfsigma$ is the Cauchy stress; $\rho$ is the mass density; $\bfv$ is the particle velocity; and $\bfb$ is the body force.

\subsection{Interface Jump Conditions and Boundary Conditions}

We use the jump operator $\llbracket \alpha \rrbracket$ to denote the difference between the limits of the discontinuous variable $\alpha$ when approaching the surface of discontinuity from either side.

The jump conditions have the form:
\begin{align}
    \label{eqn:mass-jump-main}
    \text{Mass: } & \llbracket \rho (\bfV_b - \bfv) \cdot \hat{\bfn} \rrbracket  = M
    \\
    \label{eqn:mom-jump-main}
    \text{Momentum: } & \llbracket \rho \bfv ( (\bfV_b - \bfv) \cdot \hat{\bfn} )\rrbracket + \llbracket \bfsigma \hat{\bfn} \rrbracket = M \hat{\bfv}_a
\end{align}
where $\bfV_b$ is the velocity of the surface; and $\hat\bfn$ is the surface normal.

When the ambient environment outside of the growing body at the growing surface has negligible effect on the mechanics of the body, the problem can be significantly simplified.
In particular, the domain of the problem now involves only the solid body, and the jump conditions become boundary conditions:
\begin{align}
    \label{eqn:mass-bc-main}
    \text{Mass: } & \rho (\bfV_b - \bfv) \cdot \hat{\bfn}  = M
    \quad \Longleftrightarrow \quad
    \bfV_b \cdot \hat{\bfn} = \bfv \cdot \hat{\bfn}  + \frac{M}{\rho}
    \\
    \label{eqn:mom-bc-main}
    \text{Momentum: } & \rho \bfv ((\bfV_b - \bfv) \cdot \hat{\bfn}) + 
    \bfsigma \hat{\bfn} - \bft_b = 
    M \bfv_a 
    \quad \Longleftrightarrow \quad
    \bfsigma \hat{\bfn} 
    =
    M (\bfv_a - \bfv) + \bft_b
\end{align}
where $\bft_b$ is the traction vector due to the external forces at the boundary of the growing body.

\begin{figure}[htb!]
    \includegraphics[width=\textwidth]{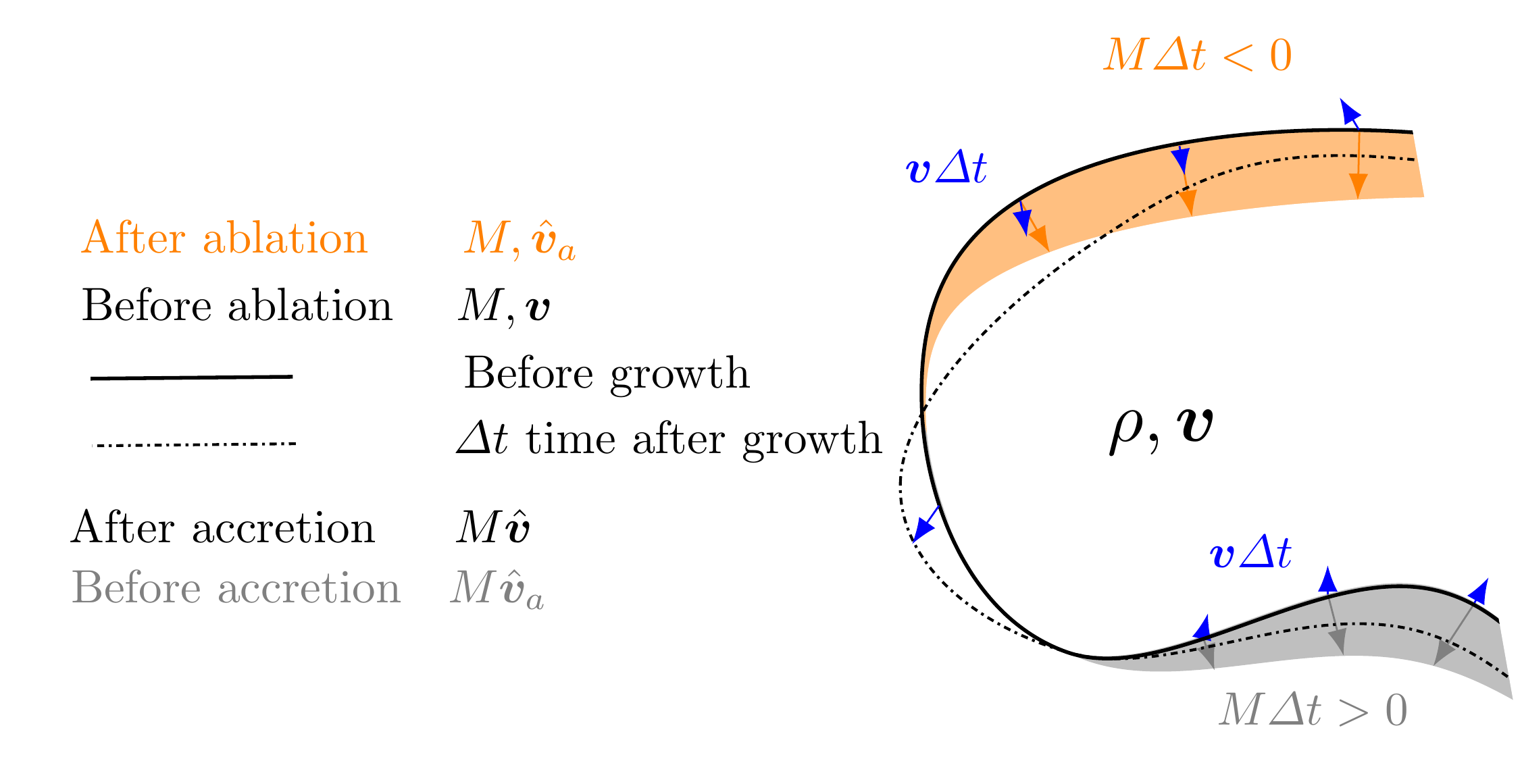}
    \caption{A schematic of the evolution of a growing body in a small interval of time $\Delta t$. The spatial location of the growing boundary depends on the growth velocity and the continuum particle velocity at the boundary.}
    \label{fig:definitions}
\end{figure}

The second form of the mass boundary condition \eqref{eqn:mass-bc-main} provides the interpretation that the boundary of body evolves due to growth as well as motion, i.e., the boundary velocity $\bfV_b$ has contributions from $M/\rho$ and $\bfv$ (Figure \ref{fig:definitions}).
The second form of the momentum boundary condition \eqref{eqn:mom-bc-main} is obtained by using \eqref{eqn:mass-bc-main}.
Interpreted physically, it states that the traction developed on the boundary is due to any imposed traction $\bft_b$ and due to the change in momentum of the added particles that occurs when their velocity transitions (instantaneously) from $\bfv_a$ to $\bfv$.

In the case of {\em slow growth}, the inertial terms $\rho \bfv ((\bfV_b - \bfv) \cdot \hat{\bfn})$ and $M \bfv_a$ are negligible relative to $\bfsigma$.
Then, the jump condition reduces to $\bfsigma \hat{\bfn} = \bft_b$, which is the standard quasistatic jump/boundary condition.

\subsection{Constitutive Response}
\label{sec:const-response}

In addition to the balance laws, we require the specification of the constitutive response of the body, corresponding here to the stress response function $\hat\bfsigma(\bfF,\nabla\bfv)$ for the initial body as well as for the added particles; in general, the stress response can be heterogeneous.

For the examples that we consider, we use an incompressible neo-Hookean stress response:
\begin{equation}
    \hat\bfsigma(\bfF) = -p\bfI + G \bfF \bfF^T
\end{equation}
where $p$ is the pressure and $G$ is the shear modulus.
In one example, we also introduce dissipation by adding the term $2\mu \sym(\nabla\bfv)$, where $\mu$ is a viscous coefficient.

\section{Kinematics and Evolution of the Elastic Deformation with Added Material}
\label{sec:kinematics}

There are 2 key issues in the kinematics of an Eulerian formulation of the growth problem.

First, we need to appropriately define the deformation state of particles at the instant of attachment.
Prior to attachment, the deformation state is irrelevant, and post-attachment deformation is governed by the standard equations of continuum mechanics.
In this regard, we highlight that the stress state of the added material is a well-defined and physically-meaningful quantity, whereas the deformation can be chosen for convenience by appropriately defining the reference at the time of attachment.

Second, the stress response of the material generally requires knowledge of the velocity gradient $\bfL:=\nabla \bfv \equiv \parderiv{\bfv}{\bfx}$ and the deformation gradient $\bfF$.
Computing $\bfL$ is straightforward, simply by applying the spatial gradient operator to the velocity field. 
The deformation gradient can be written $\bfF = \parderiv{\bfchi (\bfX , t)}{\bfX} = \left( \parderiv{\bfchi^{-1} (\bfx , t)}{\bfx} \right)^{-1}$; we notice that it requires knowledge of the motion or its inverse.
This knowledge is neither readily available nor straightforward to obtain when working in the Eulerian setting.
This problem arises in various other contexts, notably fluid-solid interaction problems (e.g. \cite{kamrin2009eulerian,kamrin2012reference,liu2001eulerian, dunne2006eulerian, sugiyama2011full}), and finite elasticity, plasticity, viscoelasticity, and biological remodeling (e.g. \cite{trangenstein1991higher,plohr1988conservative,levin2011eulerian,rajagopal2009class,clayton2019nonlinear,clayton2013nonlinear,clayton2020aplace}).

To deal with the second issue, we adapt the approach from \cite{liu2001eulerian} that formulates the evolution of $\bfF$ without requiring us to solve for the deformation map.
However, an important distinction between our strategy here and that used in \cite{liu2001eulerian} is that we introduce two additional kinematic descriptors: $\bfF_{relax}$ that quantifies the relaxed (i.e., zero stress) shape of the added particles; and $\bfF_e := \bfF \bfFr^{-1}$ that quantifies the elastic deformation.
This will be seen to allow us to pose the evolution entirely in terms of $\bfF_e$.
An important advantage in working $\bfF_e$ over $\bfF$ is that $\bfF$ must satisfy kinematic compatibility whereas $\bfF_e$ has no such requirement.

\subsection{Kinematics of the Added Material}

When new particles enter the body, we require some information about the state of the added particles.
This information can consist of specifying the following properties:
\begin{enumerate}
    \item the stress state of the added particles at attachment, $\bfsigma^*$
    \item the stress response function of the added particles, $\hat\bfsigma(\cdot)$, assumed for simplicity and without loss of generality to satisfy $\hat\bfsigma(\bfI) = \bf 0$
\end{enumerate}
We notice immediately that $\bfsigma^* = \hat\bfsigma(\bfF)$; assuming for simplicity that the stress response is invertible implies that $\bfF = \hat\bfsigma^{-1}(\bfsigma^*)$ is completely determined at the time of attachment.
Since $\bfsigma^*$ is specified independently of the state of the body, this could potentially violate kinematic compatibility at the time of attachment.

Therefore, we redefine the stress response function to be $\hat\bfsigma(\bfg) \mapsto \hat\bfsigma(\bfg \bfF_{relax}^{-1})$, where $\bfF_{relax}$ is the relaxed (i.e., zero stress) shape of the added particles.
This corresponds to redefining the reference of the added material. 
For an unattached particle, the reference can be changed in arbitrary ways and, in particular, it can be changed to ensure that we satisfy kinematic compatibility. 
We next introduce the elastic deformation $\bfF_e := \bfF \bfF_{relax}^{-1}$.
We have the obvious interpretation of $\bfF_e$ as the part of $\bfF$ that causes stress, i.e. $\bfsigma = \hat\bfsigma(\bfF_e)$.

\begin{remark}[Non-normal Growth through $\bfF_{relax}$ and $\bfF_e$]
    The kinematic descriptors $\bfF_e$ and $\bfF_{relax}$ play an important role in our approach to modeling of ``non-normal growth''.
    Non-normal growth was first discussed by \cite{skalak1982analytical}, and also further in \cite{sozio2019nonlinear, abi2019kinetics}.
    The approaches proposed there use a growth velocity vector, rather than the typical scalar normal velocity that is used to describe interface motion.
    The use of a velocity vector to describe the motion of a surface can be challenging to implement in the setting of numerical methods for free boundary problems, such as phase-field models, e.g. the formulation in \cite{agrawal2015dynamic,agrawal2015dynamic-2} which is based on a normal velocity.
    The use of $\bfF_e$ and $\bfF_{relax}$ provides an alternative: we can effectively achieve non-normal growth by using a normal interface velocity and an appropriate $\bfF_e$.
    In short, if we have a normal interface velocity, the ``shape'' of the newly added particles must be rectangular.
    However, as the material is added, we also relax the entire body using the equilibrium equation.
    In the absence of stress, the relaxed shape of the added material will correspond to $\bfF_{relax}$.
    We discuss this further with an example in Section \ref{sec:non-normal-growth}.
\end{remark}

\subsection{Evolution Equation for the Elastic Deformation}

We begin with writing down the evolution of $\bfF$, following \cite{liu2001eulerian} and others, that is posed in the form of a transport equation.
Starting from the definition of $\bfF$, we can write the material time derivative of $\bfF$ as:
\begin{equation*}
    \parderiv{\bfF(\bfX,t)}{t}
    =
    \parderiv{\ }{t}\parderiv{\bfchi(\bfX,t)}{\bfX}
    =
    \parderiv{\ }{\bfX}\parderiv{\bfx}{t}
    =
    \parderiv{\bfv}{\bfX}
    =
    \parderiv{\bfv}{\bfx}\parderiv{\bfx}{\bfX}
    =
    \left(\nabla\bfv\right) \bfF
    \Rightarrow
    \dot{\bfF}=\nabla\bfv \bfF
\end{equation*}
where $\dot{\bfg}$ represents the material time derivative of a quantity $\bfg$.
Writing this in Eulerian form, we have:
\begin{equation}
    \label{eqn:def-grad}
    \parderiv{\bfF(\bfx,t)}{t}+(\bfv\cdot\nabla)\bfF=\left(\nabla\bfv \right) \bfF
\end{equation}
This is a first-order linear tensorial advection equation that governs the time evolution of $\bfF$.

Further, $\bfF_{relax}$ is constant and independent of time for a given material particle. 
Therefore, its material derivative vanishes and we obtain:
\begin{equation}
    \label{eqn:evolution-F-relax}
    \parderiv{\bfF_{relax}}{t} + 
    (\bfv \cdot \nabla)\bfF_{relax}= 0 
\end{equation}

Right-multiplying \eqref{eqn:def-grad} by $\bfF_{relax}^{-1}$ and using \eqref{eqn:evolution-F-relax}, we obtain:
\begin{equation}
    \label{eqn:def-elastic-evolution}
    \parderiv{\bfF_e(\bfx,t)}{t}+(\bfv\cdot\nabla)\bfF_e=\left(\nabla\bfv \right) \bfF_e
\end{equation}
That is, the elastic deformation satisfies the same transport equation as the usual deformation gradient.
But the important advantage of \eqref{eqn:def-elastic-evolution} is that $\bfF_e$ has no kinematic compatibility requirements.
This makes the treatment of boundary conditions simple and straightforward, as opposed to grappling with issues of appropriate boundary conditions for \eqref{eqn:def-grad} that preserve kinematic compatibility.

\begin{remark}
    When $\dot{\bfF}_{relax}$ does not vanish identically and instead $\bfFr(t)$ is obtained in some way -- e.g., by coupling to the energy equation in a problem with thermal strain -- the evolution of the elastic deformation satisfies:
    \begin{equation}
        \dot{\bfF}_e = (\nabla\bfv)\bfF_e + \bfF_e \bfFr \dot{\overline{\left(\bfFr^{-1}\right)}}
    \end{equation}
    Equivalently, $\bfFr$ can be taken to be constant in time for a material particle, and the evolution of material properties can be incorporated through a time-dependent stress response function $\hat\bfsigma$. 
    
    While both alternatives are equivalent, the evolution of $\bfFr$ is better suited to a situation in which the evolution of material particles is easier to obtain, and the evolution of the stress response function is better suited to a situation in which the evolution of spatial locations is easier to obtain.
\end{remark}

\subsubsection{Boundary Conditions}

The evolution equation \eqref{eqn:def-elastic-evolution} for $\bfF_e$ is a transport equation.
Therefore, boundary conditions are required only at inflow boundaries, i.e., where the velocity $\bfv$ is inwards with respect to the \textit{moving} boundary.
From \eqref{eqn:mass-bc-main}, this corresponds to specifying $\bfF_e$ for added particles at growing boundaries.
No boundary condition on $\bfF_e$ are needed for the ablation boundaries.

At the inflow boundaries, we require that $\bfF_e$ be specified, which is equivalent to specifying the stress of the added particles under the assumption that the stress response function is invertible.

\subsection{Properties of the Evolution Equation}

The form of \eqref{eqn:def-elastic-evolution} suggests the method of characteristics.
Introducing the parametric variable $l$, we can write the solution of this PDE system through the solution of the following ODEs:
\begin{align}
    \label{eqn:our-char-t}
    &
    \deriv{t}{l} = 1, \quad t(l = 0) = 0 \qquad \Rightarrow \qquad t=l
    \\
    \label{eqn:our-char-x}
    &
    \deriv{x_i}{l} = v_i, \quad x_i(l=0) = c_i, \qquad i=1,2,3
    \\
    \label{eqn:our-char-F}
    &
    \deriv{(F_e)_{ij}}{l} = \parderiv{v_i}{x_k} (F_e)_{kj}, \quad (F_e)_{ij}\left(x_k(=c_k), l=t(=0)\right) = \text{initial or boundary condition} \qquad i,j,k=1,2,3
\end{align}
in a fixed Cartesian basis, and we have used index notation with implied summation.

The equations of the characteristic curves are defined by \eqref{eqn:our-char-t} and \eqref{eqn:our-char-x}:
\begin{equation}
    \hat{\bfx} = \left(\hat{x}_1(l,c_1,c_2,c_3) , \hat{x}_2(l,c_1,c_2,c_3) , \hat{x}_3 (l,c_1,c_2,c_3) , \hat{t}(l)\right)
\end{equation}
and \eqref{eqn:our-char-F} governs the behavior of $\bfF_e$ along the characteristic curves.

\subsubsection{Coincidence of Characteristic Curves and Pathlines to Track Boundary Location}

We see below that the characteristic curves $\hat{\bfx}$ of $\bfF_e$ coincide with the pathlines of the motion, and enables us to find the geometry of the body at every time.
For the parts of the boundary without growth, we simply follow the characteristic lines of particles located on the boundary; for the parts with growth, we use \eqref{eqn:mass-bc-main} together with the pathlines to determine the location of the growing boundaries. 

Following \cite{kamrin2012reference} and using the definition of the inverse of the motion $\bfX=\bfchi ^ {-1} (\bfx,t)$, we take the time derivative to find:
\begin{equation}
\label{eqn:Kamrin}
    \parderiv{\bfchi^{-1}}{t} + (\bfv\cdot\nabla) \bfchi^{-1} = 0
\end{equation}
Using the method of characteristics, the parametric ODE system of the above equation is:
\begin{align}
    &
    \label{eqn:kamrin-char-t}
    \deriv{t}{l} = 1,\quad t(l = 0) = 0 \qquad \Rightarrow \quad t=l
    \\
    &
    \label{eqn:kamrin-char-x}
    \deriv{x_i}{l} = v_i, \quad x_i(l=0) = c_i \qquad i=1,2,3
    \\
    &
    \label{eqn:chi-char}
    \deriv{\chi_i^{-1}}{l} = 0, \quad \chi_i^{-1}(x_j(=c_j) , l=t(=0)) = \text{initial or boundary condition} \qquad i,j=1,2,3
\end{align}
From \eqref{eqn:chi-char}, $\bfX=\bfchi^{-1}$ is constant along characteristic curves of \eqref{eqn:kamrin-char-x}, i.e., these characteristic curves are indeed pathlines of the motion. 
Further, comparing \eqref{eqn:kamrin-char-t} and \eqref{eqn:kamrin-char-x} to \eqref{eqn:our-char-t} and \eqref{eqn:our-char-x}, it is clear that the characteristic lines for both equations \eqref{eqn:Kamrin} and \eqref{eqn:def-elastic-evolution} are the same. 
Consequently, we conclude that the characteristic lines of $\bfF_e$ are pathlines of the motion.

\subsubsection{Solving for the Inverse Deformation Map vs. Solving for the Elastic Deformation}
\label{sec:kamrin-walkington-relation}

An interesting and powerful approach to working with solids in the Eulerian setting was presented in \cite{kamrin2012reference}.
In essence, they solve \eqref{eqn:Kamrin} for $\bfchi^{-1}$, and then compute $\bfF = (\nabla \bfchi^{-1})^{-1}$.
The relation between the approach presented here and the approach of \cite{kamrin2012reference} can be seen by first noticing that applying the spatial gradient operator to \eqref{eqn:Kamrin} recovers \eqref{eqn:def-grad}.

We assume a fixed Cartesian basis for simplicity, with uppercase indices are used to denote referential components and lowercase used to denote current components.
We start by taking the gradient of \eqref{eqn:Kamrin}:
\begin{equation}
\begin{split}
	& 
	\parderiv{\chi^{-1}_N}{t} + v_k \parderiv{\chi^{-1}_N}{x_k} = 0
	\Rightarrow
	\parderiv{\ }{x_m} \left(\parderiv{\chi^{-1}_N}{t} + v_k \parderiv{\chi^{-1}_N}{x_k} \right) = 0
	\\
	&
	\Rightarrow
	\parderiv{\ }{t} \left( \parderiv{\chi^{-1}_N}{x_m} \right) + 
	\parderiv{v_k}{x_m} \parderiv{\chi^{-1}_N}{x_k} +
	v_k \parderiv{\ }{x_k} \left( \parderiv{\chi^{-1}_N}{x_m} \right) = 0 
\end{split}
\end{equation}
Using the relation $F^{-1}_{Nk} = \parderiv{\chi^{-1}_N}{x_k}$ in the equation above, we obtain:
\begin{equation}
\begin{split}
	&
	\parderiv{F^{-1}_{Nm}}{t} + 
	\parderiv{v_k}{x_m} F^{-1}_{Nk}+
	v_k \parderiv{F^{-1}_{Nm}}{x_k} = 0
	\Rightarrow
	F_{iN} \left( \parderiv{F^{-1}_{Nm}}{t} + 
	\parderiv{v_k}{x_m} F^{-1}_{Nk}+
	v_k \parderiv{F^{-1}_{Nm}}{x_k} \right) = 0  
	\\
	&
	\Rightarrow
	\parderiv{\ }{t} \left( F_{iN} F_{Nm}^{-1} \right) - 
	\parderiv{F_{iN}}{t} F_{Nm}^{-1} + 
	v_k \parderiv{\ }{x_k} \left( F_{iN} F_{Nm}^{-1} \right)-
	v_k \parderiv{F_{iN}}{x_k} F_{Nm}^{-1} = - F_{iN} F_{Nk}^{-1} \parderiv{v_k}{x_m}
	\\
	&
	\Rightarrow
	\parderiv{F_{iN}}{t} F_{Nm}^{-1} + 
	v_k \parderiv{F_{iN}}{x_k} F_{Nm}^{-1} =  
	\parderiv{v_i}{x_m}
\end{split}
\end{equation}
Operating on both sides of the equation above by $F_{mJ}$, we recover \eqref{eqn:def-grad}.

One can therefore think of the method in \cite{kamrin2012reference} as solving \eqref{eqn:def-grad} in two steps: they first solve \eqref{eqn:Kamrin} for $\bfchi^{-1}$, and then compute $\bfF = (\nabla \bfchi^{-1})^{-1}$.
In standard continuum mechanics, where the reference and current configurations contain the same set of material particles, either of \cite{kamrin2012reference} and \cite{liu2001eulerian} could be used; there are no new particles entering the body, and hence only initial conditions are required.
Solving \eqref{eqn:Kamrin} in that setting appears easier than solving \eqref{eqn:def-grad}, as, in general, the former are three decoupled equations while the latter are nine coupled equations.
However, we highlight that when the current and reference configurations have different sets of the particles, it appears that applying inlet boundary conditions for $\bfchi^{-1}$ in \eqref{eqn:Kamrin} is challenging.
Roughly, it appears to require the construction of the reference configuration, which negates an important advantage of the Eulerian formulation.

\subsection{Postprocessing to Reconstruct the Reference Configuration and Deformation}

Our formulation above requires the coupled solution of the balance of mass and momentum along with the evolution of $\bfF_e$, summarized in \eqref{eqn:summary}.
The solution of \eqref{eqn:summary} will provide $\bfF_e(\bfx,t)$ and $\bfv(\bfx,t)$ for all $\bfx$ and $t$ in the domain of solution.
But it will \textit{not} provide $\bfF$ and $\bfF_{relax}$, which could be of interest.
A strategy to obtain $\bfF$ and $\bfF_{relax}$, once we have available $\bfF_e(\bfx,t)$ and $\bfv(\bfx,t)$, consists of the following steps:
\begin{enumerate}
    \item Choose an arbitrary time $t_0$ and define the configuration at this time as the reference configuration; consequently, $\bfF(\bfx,t_0) = \bfI$.
    
    \item Use \eqref{eqn:def-grad} to find $\bfF(\bfx,t)$ for all times of interest; notice that we have available $\bfv(\bfx,t)$, and consequently it is straightforward to solve \eqref{eqn:def-grad} numerically. We highlight that \eqref{eqn:def-grad} preserves the kinematic compatibility of $\bfF$, given a compatible field as initial data.

    \item Use the definition $\bfF_e := \bfF \bfF_{relax}^{-1}$ to find $\bfF_{relax}(\bfx,t)$. 
\end{enumerate}

\subsubsection{Connection to the Time-Evolving Reference Configuration in the Lagrangian Approach}

When new particles are added, the state of the particles must necessarily be prescribed in some form.
In the Lagrangian approach, this information can be provided in terms of the deformation state of the added particles.
Additionally, one needs to map the current state of the body at the time of attachment $\tau$ to a new reference configuration.
If this is done continuously in time as particles are added, this corresponds to a time-evolving reference configuration in which each particle is referenced through the position in the reference as well as the time of attachment of that particle.
Hence, the reference configuration is a higher order 4-dimensional manifold.

In the Eulerian setting, we could construct the evolving reference configuration by extrapolating the pathline of the particle added to the body at time $\tau$ to intersect with the $t=0$ plane, and then working with the constructed body at $t=0$ plane rather than current state of the body.  
This extrapolation depends on the attachment time $\tau$, and its representation as a 4-d manifold is denoted the ``placement map'' in \cite{tomassetti2016steady}.

\subsection{Stress Mismatch Between the Body and the Added Material}

In general, the stress state $\bfsigma^*$ of the added particles need not be consistent -- in the sense of balance of momentum -- with the stress state of the body at the region of attachment.
That is, $\bfsigma^* = \hat\bfsigma(\bfF_e)$ may not be consistent with the given traction boundary conditions or the stress state of the body at the growing boundary.

This is resolved by the evolution of the deformation and stress in both the added particles and the body following the balance of momentum, which is solved simultaneously with the kinematic equation in \eqref{eqn:summary}.
In the quasistatic setting, this occurs instantaneously through the relaxation, and in the dynamic setting this occurs through elastodynamics.

To illustrate this scenario, consider two brief examples.
\begin{enumerate}
    \item 
    Consider hot particles attaching to a cool body with traction free boundaries, such as in additive manufacturing.
    At the instant of attachment, the hot particles are assumed to have no stress.
    However, the hot particles cool instantaneously upon attachment, assuming that heat transfer is fast, and there is consequent thermal strain.
    
    This entire process could be coarse-grained as the attachment of cool and stressed particles -- $\bfF_e \neq \bfI$ and $\bfsigma^* \neq \bf 0$ to compensate for the thermal strain -- on a traction free substrate.

    \item
    Consider deposition on a body subject to a shear traction, such as in biological tissue next to a flowing fluid.
    The particles in the fluid are typically stress-free, but undergo shear as soon as they attach to body.
    While the details of the attachment and subsequent shearing  could be a complex time-dependent process, we can coarse-grain this process to apply the deformation and shear stress on the added particles at the instant of attachment.

\end{enumerate}

While the quasistatic setting is typically more appropriate for growth, the instantaneous relaxation of the stress can lead to discontinuities in the motion.
While these discontinuities are physically reasonable under the assumption of instantaneous relaxation, they make it challenging to find solutions.
In the interest of simplicity, we follow the strategy of \cite{kamrin2012reference} and add a viscous term ($2 \mu \sym(\nabla \bfv)$) to the stress response, as in Section \ref{sec:const-response},
to regularize the motion and make it continuous.
The quasistatic solution is defined as the limit of $\mu \to 0$.


\section{Summary of the Proposed Approach}
\label{sec:summary}

Our approach requires the solution of the balances of mass and momentum, and the transport of the elastic deformation.
In summary, we solve the 3 coupled equations given by:
\begin{equation}
\label{eqn:summary}
    \deriv{\ }{t} \begin{Bmatrix}
     \rho \\ \bfv \\ \bfF_e 
    \end{Bmatrix}
    =
    \begin{Bmatrix}
    -\rho \divergence\bfv \\ \frac{1}{\rho}\divergence \hat\bfsigma(\bfF_e) \\ (\nabla\bfv) \bfF_e 
    \end{Bmatrix}
\end{equation}
where $\deriv{\ }{t}$ is the material time derivative.

\section{Example: Non-normal Growth Through the Elastic Deformation}
\label{sec:non-normal-growth}

\begin{figure}[ht!]
    \includegraphics[width=\textwidth]{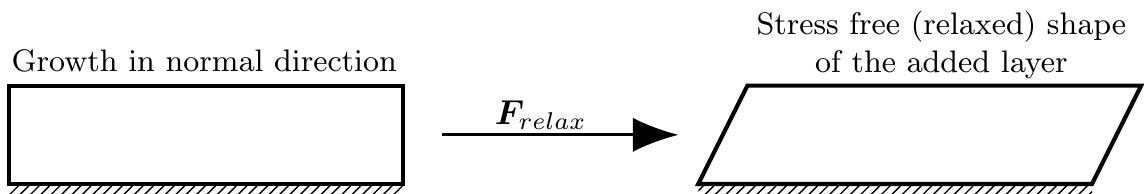}
    \caption{Modeling non-normal growth using $\bfF_e$ and  $\bfF_{relax}$.}
    \label{fig:tangential_growth}
\end{figure}

We consider the modeling of non-normal growth use the kinematic descriptors $\bfFr$ and $\bfF_e$.
\cite{skalak1982analytical} discussed the case of non-normal growth, shown in Figure \ref{fig:tangential_growth} in the right panel, wherein 
the material grows in an inclined fashion.

Assume that the surface is traction free.
If the growth velocity is normal to the surface and the layer is stress free, the ``shape'' of the added layer will be rectangular, and remain rectangular as it is in equilibrium.
However, if we set $\bfF_e = \begin{bmatrix}
    1 & -\alpha \\
    0 & 1 \\
\end{bmatrix}$ at the growing boundary, the added layer is not stress free; it will consequently relax to a sheared shape when we solve for equilibrium with zero traction at the surface, as in Figure \ref{fig:tangential_growth}.

\subsection{Reduced set of the governing equations}

Consider the growth of a neo-Hookean material with normal growth velocity $V_G=\frac{M}{\rho}$ clamped on a substrate. 
We make the following simplifying assumptions:
\begin{itemize}
	\item The material is incompressible. 
	Further, based on the definition of the problem, we assume that the velocity field has the form $\bfv = v_1 (x_1,x_2,t)\hat \bfe_1$ ($v_2(x_1,x_2,t)=0$). 
	Then, the bulk mass balance \eqref{eqn:bulk-balance-mass} reduces to \eqref{eqn:mass-non-normal}.
	
	\item 
	One consequence of the previous assumption is that the height of the body is $H=H(t)$ and the top growing boundary remains flat with normal $\hat e_2$, so the velocity of the growing top boundary is \eqref{eqn:growing-vel-bc-non-normal}. 
	
	\item The growth is slow, and inertial terms are neglected. Also there are no body forces. This reduces the balance of momentum \eqref{eqn:bulk-balance-momentum} and \eqref{eqn:mom-bc-main} to \eqref{eqn:mom-non-normal} and \eqref{eqn:bc-growing-mom-non-normal}, respectively. 
	
	\item The constitutive law for the neo-Hookean material with addition of a viscous term in the constitutive law to make the motion continuous is \eqref{eqn:const-law-non-normal}, where the pressure $p$ is an unknown Lagrange multiplier function. 
	
	\item All of the added and original material is identical, and hence the referential density of the added particles is equal to the constant referential density of the body. Combined with the incompressibility assumption, we have the current density is also constant.
	
	\item We ignore any thermal effects.
	
	\item To avoid discontinuity in the motion due to the slow growth assumption, we add viscous term $2 \mu \sym(\nabla \bfv)$ to the constitutive law of the material. 
	
\end{itemize}

The model equations with the boundary conditions is given by:
\begin{align}
    &
    \label{eqn:mass-non-normal}
    \divergence(\bfv)
    = 
    \parderiv{v_1}{x_1} 
    =
    0 
    \quad \Rightarrow \quad
    v_1 = v_1(x_2,t)
    \quad \text{on} \quad
    0 < x_2 < H(t)
    \\
    &
    \label{eqn:fixed-wall-bc-non-normal}
    \bfV_b = \bfv = 
    v_1 \hat\bfe_1 + v_2 \hat\bfe_2 = 
    \mathbf{0}
    \quad \Rightarrow \quad
    v_1 = v_2 = 0
    \quad \text{at} \quad
    x_2=0
    \\
    &
    \label{eqn:growing-vel-bc-non-normal}
    \bfV_b \cdot \hat\bfn = \bfv \cdot \hat\bfn + V_G
    \quad \text{at} \quad
    x_2 = H(t)
    \ , \
    \hat\bfn = \hat\bfe_2 
    \quad \Rightarrow \quad
    V_{b_2} = V_G
    \quad \Rightarrow \quad
    H(t) = V_G t
    \\
    &
    \label{eqn:mom-non-normal}
    \divergence(\bfsigma) = 
    \left( \parderiv{\sigma_{11}}{x_1} + \parderiv{\sigma_{12}}{x_2} \right) \hat\bfe_1 +
    \left( \parderiv{\sigma_{12}}{x_1}+ \parderiv{\sigma_{22}}{x_2} \right) \hat\bfe_2
    = \mathbf{0}
    \quad \text{on} \quad
    0 < x_2 < H(t) 
    \\
    &
    \label{eqn:bc-growing-mom-non-normal}
    \bfsigma \hat\bfn = 
    \sigma_{12} \hat\bfe_1 + \sigma_{22} \hat\bfe_2 = 
    \mathbf{0} 
    \quad \text{at} \quad
    x_2=H(t)
    \\
    &
    \label{eqn:Fdot_non-normal}
    \parderiv{\ }{t}
    \begin{bmatrix}
    F_{e_{11}} & F_{e_{12}} \\
    F_{e_{21}} & F_{e_{22}} \\
    \end{bmatrix}
    =
    \begin{bmatrix}
    0 & \parderiv{v_1}{x_2} \\
    0 & 0 \\
    \end{bmatrix}
    \begin{bmatrix}
    F_{e_{11}} & F_{e_{12}} \\
    F_{e_{21}} & F_{e_{22}} \\
    \end{bmatrix}
    =
    \begin{bmatrix}
    \parderiv{v_1}{x_2} F_{e_{21}} & \parderiv{v_1}{x_2} F_{e_{22}} \\
    0 & 0
    \end{bmatrix}
    \quad \text{on} \quad
    0 < x_2 < H(t)
    \\
    &
    \label{eqn:BC-F-non-normal}
    \bfF_e = 
    \begin{bmatrix}
        1 & -\alpha \\
        0 & 1 \\
    \end{bmatrix}
    \quad \text{at} \quad
    t=t_{att} = \frac{x_2}{V_G}
    \\
    &
    \label{eqn:const-law-non-normal}
    \hat\bfsigma(\bfF_e,\nabla\bfv) =
    -p \bfI + G \bfF_e \bfF_e^T + 2 \mu \sym(\nabla \bfv)
\end{align}

\subsection{Solution }

The evolution equation for $F_{e_{21}}$ is homogeneous and the boundary condition for it is equal to zero, so $F_{e_{21}}(x_1,x_2,t) = 0$. 
Also, the evolution equation for the diagonal components of $\bfF_e$ will be also homogeneous and the boundary condition on them is constant in time and $x_1$, so all the characteristic curves transfer the same value and keep them constant as the equations are homogeneous.
So $F_{e_{11}}(x_1,x_2,t) = F_{e_{22}}(x_1,x_2,t) = 1$ which is the value at the boundary.
Therefore, $\bfF_e$ has the form:
$\bfF_e(x_1,x_2,t) = \begin{bmatrix}
    1 & F_{e_{12}}(x_1,x_2,t) \\
    0 & 1 \\
    \end{bmatrix}$.
Moreover, we assume that $F_{e_{12}} = F_{e_{12}}(x_2,t)$ and it does not depend on $x_1$. Substituting $\bfF_e(x_1,x_2,t) = \begin{bmatrix}
    1 & F_{e_{12}}(x_2,t) \\
    0 & 1 \\
    \end{bmatrix}$ 
in the above set of the equations, the system of the equations for computing $v_1(x_2,t),\ F_{e_{12}}(x_2,t),\ p(x_1,x_2,t)$ reduces to the following form:
\begin{align}
    &
    \parderiv{F_{e_{12}}}{t} =\parderiv{v_1}{x_2}
    , \quad 
    \text{ with } F_{e_{12}} = -\alpha \quad \text{at} \quad t=t_{att} = \frac{x_2}{V_G}
    \\
    &
    -\parderiv{p}{x_1} + G \parderiv{F_{e_{12}}}{x_2} + \mu \parderivsec{v_1}{x_2} = 0
    , \quad
    \parderiv{p}{x_2} = 0
    , \quad
    \text{ with } \sigma_{12} = \sigma_{22} = 0 \quad 
    \text{at} \quad x_2 = H(t) = V_Gt 
\end{align}
The solution is:
\begin{align}
    &
    v_1(x_1,x_2,t) = V_G \alpha 
    \left( e^{-\frac{G}{\mu}\left(t-\frac{x_2}{V_G}\right)} - 
    e^{ -\frac{G}{\mu} t }
    \right) 
    \xlongequal{\mu \to 0}
    0
    \\
    &
    F_{e_{12}}(x_1,x_2,t)=
    - \alpha  \exp{\left(-\frac{G}{\mu}\left(t-\frac{x_2}{V_G}\right)\right)}
    \xlongequal{\mu \to 0}
    0
    \\
    &
    p(x_1,x_2,t)=G \quad \Rightarrow \quad 
    \sigma_{12} (x_1,x_2,t) = \sigma_{22} (x_1,x_2,t) = 0
\end{align}
The quasi-static solution can be computed in the limit of $\mu \to 0$. 
In this case, $F_{e_{12}}$ is equal to zero everywhere and for all times after attachment.
This implies that the solution of the equilibrium equation immediately after attachment causes $F_{e_{12}}$ to change instantaneously from $-\alpha$ to zero.
That is, the particles shear over immediately after attachment, driven by a combination of a choice of $\bfF_e$ and momentum balance.


\section{Example: Compatible Deformations for Growth with Thermal Strain Variation}
\label{sec:compatible-deformations}

\begin{figure}[ht!]
    \centering
    \includegraphics[width=\textwidth]{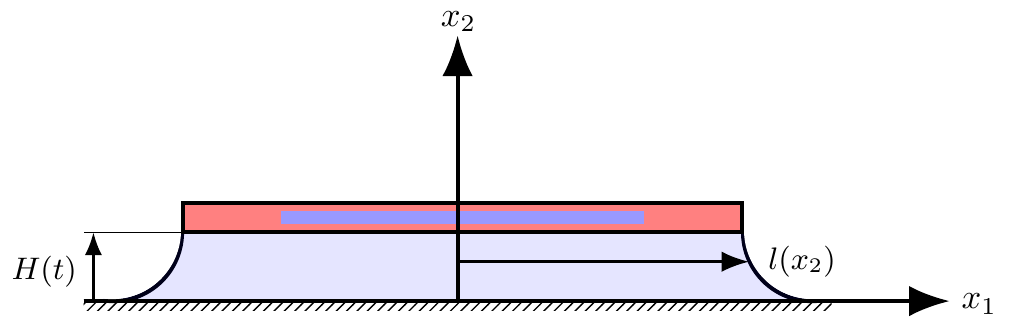}
    \caption{Modeling thermal growth. The light blue domain is the existing body that is clamped to a substrate; the red domain is the added material; and the dark blue domain is the stress-free shape of the added material after cooling.}
    \label{fig:thermal_growth}
\end{figure}

We consider growth in a setting that roughly models the process of additive manufacturing with large thermal gradients.
This example illustrates the compatibility of deformations in a setting with thermal strain.
Growth occurs by the deposition of a stress-free layer of material that is hot with temperature $T_h$ on a cool body with temperature $T_c$; these are shown schematically by the red and light blue domains in Figure \ref{fig:thermal_growth}.
As material is added to the body, it cools and shrinks, but is clamped to the substrate at the base.

In the limit of infinite thermal conduction, the added layer instantly cools to $T_c$ upon deposition, and consequently induces thermal strain.
That is, the added material has a relaxed shape that would, typically, correspond to a contraction.
This thermal strain value depends on $T_c$ and $T_h$.

We start by assuming that the added material is stress-free at the instant of attachment. 
However, immediately after attachment, the added material instantly cools; in the absence of constraints due to the existing body, it would develop a thermal strain (shown by the dark blue rectangle in Figure \ref{fig:thermal_growth}).
However, the constraints on the added material due to the attachment to the existing body induce an elastic deformation, assumed for simplicity to have the form $\bfF_e = \frac{1}{\alpha} \bfI$, where $\alpha$ is a function of $T_h$ and $T_c$.

Obviously, the added material is not in equilibrium with the body and traction free condition at the growing boundary. 
Under the quasistatic assumption, the deformation of the newly added layer occurs instantly after attachment, and can be computed by solving the equilibrium equation for the added layer and the body.
However, we can also regularize the problem by using a finite thermal conductivity.

The assumptions that we make are as follows:
\begin{itemize}
    \item  
    We assume a flat traction-free surface and neglect any nonuniformity at the ends.
    
    \item 
    The material is incompressible. Then, the bulk mass balance reduces to \eqref{eqn:mass-thermal}, and the constitutive law for a hyperelastic material is \eqref{eqn:const-law-thermal}, where the pressure $p$ is an unknown Lagrange multiplier function. 

    \item 
    The growth is slow, and inertial terms are neglected. Also there are no body forces.  This reduces the balance of momentum to \eqref{eqn:mom-thermal} and \eqref{eqn:bc-growing-mom-thermal}. 
\end{itemize}

So the governing equations are:
\begin{align}
    &
    \label{eqn:mass-thermal}
    \divergence(\bfv) = 0 
    \quad \text{on} \quad
    0 < x_2 < H(x_1,t)
    \\
    &
    \label{eqn:fixed-wall-bc-thermal}
    \bfV_b = \bfv = 
    v_1 \hat\bfe_1 + v_2 \hat\bfe_2 = 
    \mathbf{0}
    \quad \Rightarrow \quad
    v_1 = v_2 = 0
    \quad \text{at} \quad
    x_2=0
    \\
    &
    \label{eqn:growing-vel-bc-thermal}
    \bfV_b \cdot \hat\bfn = \bfv \cdot \hat\bfn + V_G
    \quad \text{at} \quad
    x_2 = H(t)
    \ , \
    \hat\bfn = \hat\bfe_2 
    \quad \Rightarrow \quad
    V_{b_2} = V_G+v_2
    \\
    &
    \label{eqn:mom-thermal}
    \divergence(\bfsigma) = 
    \left( \parderiv{\sigma_{11}}{x_1} + \parderiv{\sigma_{12}}{x_2} \right) \hat\bfe_1 +
    \left( \parderiv{\sigma_{12}}{x_1}+ \parderiv{\sigma_{22}}{x_2} \right) \hat\bfe_2
    = \mathbf{0}
    \quad \text{on} \quad
    0 < x_2 < H(t) 
    \\
    &
    \label{eqn:bc-growing-mom-thermal}
    \bfsigma \hat\bfn = 
    \sigma_{12} \hat\bfe_1 + \sigma_{22} \hat\bfe_2 = 
    \mathbf{0} 
    \quad \text{at} \quad
    x_2=H(t)
    , \quad 
    \bfsigma \hat\bfn =
    \mathbf{0} 
    \quad \text{at} \quad
    x_1=\pm l(x_2), \hat\bfn=\hat\bfn(x_2)
    \\
    &
    \label{eqn:Fdot_thermal}
    \parderiv{\bfF_e}{t} + 
    v_1 \parderiv{\bfF_e}{x_1} + 
    v_2 \parderiv{\bfF_e}{x_2}
    =
    \nabla \bfv \bfF_e
    \quad \text{on} \quad
    0 < x_2 < H(t)
    \\
    &
    \label{eqn:BC-F-thermal}
    \bfF_e = \frac{1}{\alpha}\bfI 
    \quad \text{at} \quad
    t=t_{att} = \frac{x_2}{V_G}
    \\
    &
    \label{eqn:const-law-thermal}
    \hat\bfsigma(\bfF_e) =
    -p \bfI + G (\bfF_e) (\bfF_e)^T,
\end{align}

In general, the solution for the above set of the equations depends on all independent variables $x_1, x_2$ and $t$.
Even with several reasonable simplifying assumptions, we are unable to solve these in closed-form, and hence we leave the (seemingly straughtforward) numerical solution for future work.

\section{Example: Growth due to Accretion with Inertia and Shear}
\label{sec:shear}

We consider growth that occurs on a traction-free surface, but with momentum transfer due to the interaction between moving added material and a stationary body.
This serves as a first approximation to a number of interesting situations.
For instance, this could model the growth of biological tissue while subject to shear and normal tractions from fluid motion over the growing surface.
Alternatively, it could model the deposition of a sticky material on the growing boundary of a body, wherein the relative velocity between the sticky material and the body causes a shear force on the body due to the change in the momentum of the added material.

\begin{figure}[htb!]
    \centering
    \includegraphics[width=\textwidth]{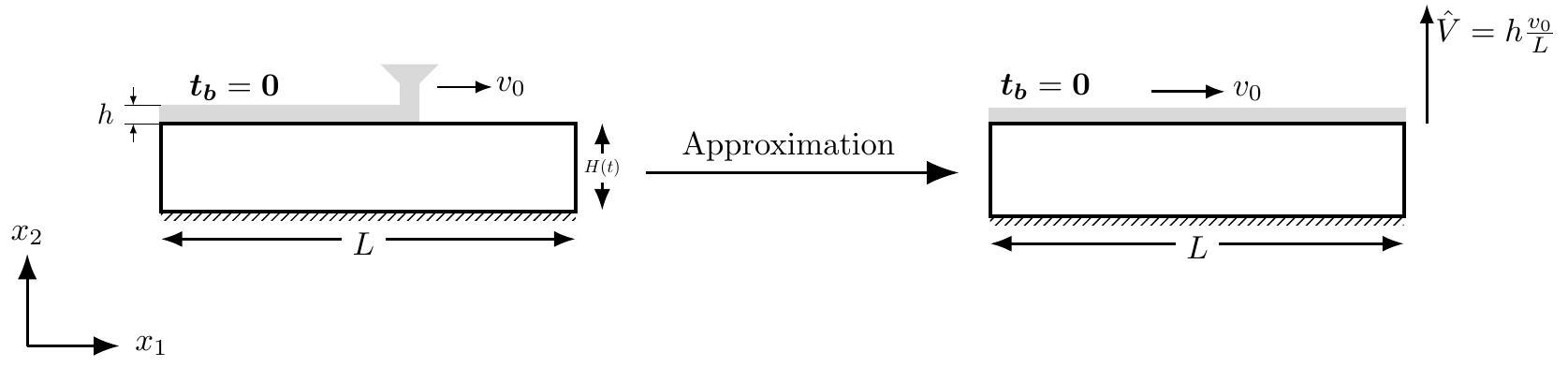}
    \caption{Schematic view of a realistic and idealized Fused Deposition Modeling (FDM) process.}
    \label{fig:FDM}
\end{figure}

Specifically, we consider here the additive manufacturing process of fused deposition modeling.
A sticky material is deposited in layers over an existing body, with each layer of thickness $h$ related to the nozzle radius.
The nozzle moves with the velocity $v_0$ to the right (figure \ref{fig:FDM}). 
The particles are idealized to have a nozzle exit velocity -- relative to the lab frame -- of $v_0 \hat\bfe_1$ with no vertical component.

To model this as a continuous growth process, we coarse-grain the details as follows.
We ignore the process of horizontal motion of the nozzle in defining the rate of mass addition, and assume a uniform (in space and time) mass rate of growth per unit area given by $M = \rho h \frac{v_0}{L}$.
However, the horizontal velocity of the added material is important in defining the shear traction that acts on the body, and will be accounted for through the momentum jump condition \eqref{eqn:mom-bc-main}.

\subsection{Formulation and Reduced Governing Equations}

Consider a long rectangular body of a neo-Hookean material in 2 dimensions with initial height $H_0$ ($H(t) \ll L$), which is clamped at $x_2 = 0$, see \eqref{eqn:fixed-wall-bc-shear}, and initially stress-free, see \eqref{eqn:Fe-bc-Ic-shear}. 
Initially stress-free particles, see \eqref{eqn:Fe-bc-Ic-shear}, with an initial velocity of $\hat{\bfv}_a = v_0 \hat\bfe_1$ are attached to the body at the upper growing boundary.
The external traction $\bft_b$ at this boundary is zero. 

The following assumptions are considered for the formulation of the problem:
\begin{itemize}
	\item The body is long in the $x_1$ direction, so we neglect the dependence on $x_1$, and approximate all quantities to be only functions of $x_2$ and $t$. Consequently, all derivatives with respect to $x_1$ are $0$.
	
	\item The material is incompressible; together with the previous assumption, the bulk mass balance \eqref{eqn:bulk-balance-mass} reduces to \eqref{eqn:mass-shear}.
	
	\item The growth is slow, and inertial terms the momentum balance are neglected: i.e., $\rho \Deriv{\bfv}{t} \ll \divergence(\bfsigma)$ and $M \bfv \ll \bfsigma \hat\bfn$. Also there are no body forces. This reduces the balance of momentum in bulk \eqref{eqn:bulk-balance-momentum} to \eqref{eqn:mom-shear}.
	
	\item The external traction at the growing boundary is zero.
	However, the added material has a non-negligible initial velocity of $\hat{\bfv}_a = v_0 \hat\bfe_1$.
	Strictly, the momentum transfer at the growing boundary is $M (\hat{\bfv}_a - \bfv)$, but we assume that $|\bfv| \ll  |\hat{\bfv}_a|$, so the boundary condition of the momentum balance \eqref{eqn:mom-bc-main} reduces to \eqref{eqn:bc-growing-mom-shear}. 
	
	\item The constitutive law for the neo-Hookean material is \eqref{eqn:const-law-shear}, where the pressure $p$ is an unknown Lagrange multiplier function. 
	
	\item All the material considered is the same. Further, we assume that conditions are isothermal and the material is incompressible.
	Consequently, the referential density of the added particles is equal to the constant referential density of the body, and the current density is also constant.
	
\end{itemize}

So, the set of the governing equations is:
\begin{align}
    &
    \label{eqn:mass-shear}
    \divergence(\bfv)
    =
    \parderiv{v_2}{x_2}
    = 
    0 
    \quad \Rightarrow \quad
    v_2 = v_2(t)
    \quad \text{on} \quad
    0 < x_2 < H(t)
    \\
    &
    \label{eqn:fixed-wall-bc-shear}
    \bfV_b = \bfv = 
    v_1 \hat\bfe_1 + v_2 \hat\bfe_2 = 
    \mathbf{0}
    \quad \Rightarrow \quad
    v_1 = v_2 = 0
    \quad \text{at} \quad
    x_2=0
    \quad \Rightarrow \quad
    v_2(t)
    = 0
    \\
    &
    \label{eqn:growing-vel-bc-shear}
    \bfV_b \cdot \hat\bfn = \bfv \cdot \hat\bfn + \frac{M}{\rho}
    \quad \text{at} \quad
    x_2 = H(t)
    \ , \
    \hat\bfn = \hat\bfe_2 
    \quad \Rightarrow \quad
    V_{b_2} = h \frac{v_0}{L}
    \\
    &
    \label{eqn:mom-shear}
    \divergence(\bfsigma) = 
    \left( \parderiv{\sigma_{12}}{x_2} \right) \hat\bfe_1 +
    \left( \parderiv{\sigma_{22}}{x_2} \right) \hat\bfe_2
    = \mathbf{0}
    \quad \text{on} \quad
    0 < x_2 < H(t) 
    \\
    &
    \label{eqn:bc-growing-mom-shear}
    \bfsigma \hat\bfn = 
    \sigma_{12} \hat\bfe_1 + \sigma_{22} \hat\bfe_2 = 
    \bft_b + M \hat{\bfv}_a = 
    \rho h \frac{v_0^2}{L} \hat\bfe_1 
    \quad \text{at} \quad
    x_2=H(t)
    \\
    &
    \label{eqn:Fdot_shear}
    \parderiv{\ }{t}
    \begin{bmatrix}
    F_{e_{11}} & F_{e_{12}} \\
    F_{e_{21}} & F_{e_{22}} \\
    \end{bmatrix}
    =
    \begin{bmatrix}
    0 & \parderiv{v_1}{x_2} \\
    0 & 0 \\
    \end{bmatrix}
    \begin{bmatrix}
    F_{e_{11}} & F_{e_{12}} \\
    F_{e_{21}} & F_{e_{22}} \\
    \end{bmatrix}
    =
    \begin{bmatrix}
    \parderiv{v_1}{x_2} F_{e_{21}} & \parderiv{v_1}{x_2} F_{e_{22}} \\
    0 & 0
    \end{bmatrix}
    \quad \text{on} \quad
    0 < x_2 < H(t)
    \\
    &
    \label{eqn:Fe-bc-Ic-shear}
    \bfF_e (t=0, x_2<H_0) = \bfI,
    \quad 
    \bfF_e (t=t_{att}, x_2) = \bfI
    \\
    &
    \label{eqn:const-law-shear}
    \hat\bfsigma(\bfF_e) =
               -p \bfI + G \bfF_e \bfF_e^T 
\end{align}

\subsection{Solution}

From the continuity equation \eqref{eqn:mass-shear}, the boundary condition \eqref{eqn:fixed-wall-bc-shear}, and the continuity of velocity on the interface between the initially existing body and the accreted part of the body, we see that $v_2 \equiv 0$ throughout the growing body.
Then, the top surface of the body is located at $x_2 = H=H(t)$, and the top growing boundary remains horizontal straight line with $\hat\bfn=\hat e_2$.
So the upper growing boundary can be located by integrating \eqref{eqn:growing-vel-bc-shear}, and is at ${x_b}_2 (t) = H(t) = H_0 + \frac{h v_0}{L} t$.
Consequently, the time of attachment of the newly added particles is $t_{att} = \frac{x_2 - H_0}{hv_0/L}$.

In the above set of the equations, the initial body is at rest, and the added material is also assumed to be stress free, so $\bfF_e(t=0) = \bfF_e (t=t_{att}) = \bfI$. 
But the added material imposes a uniform shear traction on the upper growing surface, which is not consistent with the stress state of the body.
Using the quasistatic assumption (i.e., neglecting the bulk inertial effects in the balance of momentum), the initial body and the added particles both respond instantaneously to the traction at the boundary and there will be a jump (in time) in $\bfF_e$ at $t=0$ for $x_2<H_0$ and $t=t_{att}$ for $x_2\geq H_0$. 
The discontinuity is readily handled by solving for equilibrium, and the value of $\bfF_e$ immediately after deposition begins is given by:
\begin{equation}
    \bfF_e = \begin{bmatrix}
    1 & \frac{M v_0}{G} \\
    0 & 1
    \end{bmatrix}
    \quad \text{at} \quad
    t=0^+, \ 0<x_2<H_0
\end{equation}
Then using \eqref{eqn:Fdot_shear} together with the condition above for the region $0<x_2<H_0$, the deformation gradient in the initially existing body is:
    \begin{equation}
    \bfF_e (x_2 , t) =
        \begin{bmatrix}
        1 & F_{e_{12}} \\
        0 & 1
        \end{bmatrix} 
        \qquad \text{in} \qquad 0 < x_2 < H_0
    \end{equation}
where $F_{e_{12}}$ satisfies the condition $\parderiv{F_{e_{12}}}{t} = \parderiv{v_1}{x_2}$, and the corresponding stress state is:
\begin{equation}
\label{eqn:sigma_shear_r1}
    \bfsigma (x_2 , t)=
    \begin{bmatrix}
    -p + G (1+F_{e_{12}}^2) & G F_{e_{12}} \\
    G F_{e_{12}} & -p+G
    \end{bmatrix} 
    \qquad \text{in} \qquad 0 < x_2 < H_0
\end{equation}

In addition, the boundary condition for $\bfF_e$ at $t=t_{att}$ in the growing part of the body ($x_2 \geq H_0$) has to satisfy the following conditions:
\begin{itemize}
    \item The incompressibility of the material,
    
    \item The traction boundary condition  $\bft = M v_0 \hat\bfe_1$ at the upper layer with $\hat\bfn=\hat\bfe_2$, 
    
    \item The modeling assumption that the deformation of the added material at the time of attachment is such that there is no deformation in the tangential plane 
    so $\bfF_e \bfl =\bfl$ for any $\bfl$ satisfying $\bfl\cdot\hat\bfn = 0$.
\end{itemize}
Then, the boundary condition for $\bfF_e$ at the growing boundary is chosen to be:
\begin{equation*}
    \bfF_e = \begin{bmatrix}
    1 & \frac{M v_0}{G} \\
    0 & 1
    \end{bmatrix}
    \quad \text{at} \quad
    x_2 = H(t)  
\end{equation*}
which is compatible with all these assumptions. 

Using the above boundary and initial conditions for $\bfF_e$, together with the governing equations, and the continuity of the velocity and the normal component of the stress at the interface between the initially existing body and the added materials, the solution has the form:
\begin{align*}
    &
    \bfsigma (x_2 ,t) = 
    \begin{bmatrix}
    \frac{(M v_0)^2}{G} 
    & 
    M v_0  
    \\
    M v_0 
    & 
    0
    \end{bmatrix}
    , \quad
    \bfF_e (x_2 ,t) = 
    \begin{bmatrix}
    & 1 & \frac{M v_0}{G} \\
    & 0 & 1
    \end{bmatrix}
    \quad , \quad
    v_1(x_2 , t) =0
    \quad \text{on} \quad
    0 < x_2 < H(t)
\end{align*}
It is easy to check this by substituting in \eqref{eqn:mom-shear} -- \eqref{eqn:const-law-shear}.

We notice, as expected, that the body deforms uniformly with a shear stress.


\section{Discussion}
\label{sec:conclusion}

In this paper, we model the surface growth of solid bodies using an Eulerian approach.
An important advantage of the Eulerian approach is that we do not need to explicitly compute the time-evolving reference configuration.
On the other hand, the solid stress response requires a knowledge of the deformation gradient $\bfF$, which is challenging to obtain in the Eulerian setting.
We address this issue by introducing the additional kinematic descriptors of the zero stress deformation $\bfFr$ and the elastic deformation $\bfF_e$.
We then eliminate $\bfF$ and $\bfFr$ and formulate a model whose only kinematic descriptor is $\bfF_e$.
This formulation provides some important advantages: 
(1) transparent specification of boundary conditions for $\bfF_e$ directly from the stress state of the added particles;
(2) boundary conditions for $\bfF_e$ come \textit{solely} from the stress state of the added particles, with no restrictions due to kinematic compatibility;
and (3) the ability to model complex phenomena such as non-normal growth using only a normal interface velocity, making them amenable to standard free boundary methods.

The evolution of $\bfF_e$ has the structure of a transport equation, providing a natural role for boundary conditions to describe the state of added material at inflow boundaries and no boundary conditions are imposed at outflow boundaries.

A primary motivation for formulating the Eulerian approach is to enable future work on numerical solutions, following the Lagrangian approach pioneered in \cite{von2020morphogenesis}.
A potential advantage of Eulerian methods that have been developed in the fluid-structure interaction (FSI) literature, e.g., \cite{kamrin2009eulerian,kamrin2012reference,liu2001eulerian, dunne2006eulerian, sugiyama2011full,kamensky2017immersogeometric}, is that we can solve such problems on a fixed mesh for an evolving body, rather than an evolving reference configuration and deforming mesh that needs to be updated or re-meshed frequently.
A numerical scheme would enable the careful study of various interesting systems identified recently, e.g. \cite{zurlo2017printing, zurlo2018inelastic, truskinovsky2019nonlinear} that focus on the problem of near-net-shape manufacturing in the presence of incompatibility and residual stress, as well as \cite{swain2018biological} that discusses incompatibility in both surface and bulk growth.


\begin{acknowledgments}
    We thank Tal Cohen for introducing us to this topic and useful discussions; Rohan Abeyaratne, Timothy Breitzman, and Tony Rollett for useful discussions; the anonymous reviewers for suggestions that improved the paper; and NSF (CMMI MOMS 1635407, DMREF 1628994, DMS 1729478, and DMS 2012259), ARO (MURI W911NF-19-1-0245), ONR (N00014-18-1-2528), and BSF (2018183) for financial support.
\end{acknowledgments}


\bibliographystyle{alpha}
\bibliography{growth-refs}

\end{document}

%% file: preamble.tex
\usepackage{isomath}
\usepackage{amsmath,amsthm}
\usepackage{amsbsy}
\usepackage{amssymb}
\usepackage{amscd}
\usepackage{amsfonts}
\usepackage{stmaryrd}
\usepackage{siunitx}
\usepackage{euscript}
\usepackage[utf8]{inputenc}
\usepackage[T1]{fontenc}
\usepackage{newtxtext} 
\everymath{\displaystyle}
\usepackage{exscale}

\usepackage{graphicx}
\usepackage{boxedminipage}
\usepackage{calc}
\usepackage[usenames,dvipsnames]{xcolor}
\usepackage[caption=false]{subfig}

\usepackage{setspace}
\usepackage{enumitem}
\setitemize{noitemsep,topsep=0pt,parsep=0pt,partopsep=0pt}
\setenumerate{noitemsep,topsep=0pt,parsep=0pt,partopsep=0pt}
\setdescription{noitemsep,topsep=0pt,parsep=0pt,partopsep=0pt}

\usepackage[,colorlinks,urlcolor=blue,citecolor=black]{hyperref}
\usepackage[normalem]{ulem}

\usepackage[small]{titlesec}

\titlespacing*{\section}{0pt}{12pt plus 4pt minus 2pt}{2pt plus 2pt minus 2pt}
\titlespacing*{\subsection}{0pt}{12pt plus 4pt minus 2pt}{2pt plus 2pt minus 2pt}
\titlespacing*\subsubsection{0pt}{12pt plus 4pt minus 2pt}{2pt plus 2pt minus 2pt}
\titlespacing*\paragraph{0pt}{12pt plus 4pt minus 2pt}{2pt plus 2pt minus 2pt}

\makeatletter

    \renewcommand*{\p@subsection}{}
    
    \renewcommand*{\p@subsubsection}{}
\makeatother

\input{definitions}

%% file: definitions.tex
\usepackage{isomath}
\usepackage{amsmath}
\usepackage{amssymb}
\usepackage{amscd}
\usepackage{amsfonts}
\usepackage{extarrows}

\newtheorem{remark}{Remark}[section]

\newcommand{\bfchi}{\mathbold {\chi}}

\newcommand{\bfsigma}{\mathbold {\sigma}}

\DeclareMathOperator{\divergence}{div}

\DeclareMathOperator{\sym}{sym}

\newcommand{\parderiv}[2]{\frac{\partial #1}{\partial #2}}
\newcommand{\parderivsec}[2]{\frac{\partial^2 #1}{\partial #2 ^2}}
\newcommand{\dm}{\ \mathrm{d}}
\newcommand{\Dm}{\ \mathrm{D}}
\newcommand{\deriv}[2]{\frac{\dm #1}{\dm #2}}
\newcommand{\Deriv}[2]{\frac{\Dm #1}{\Dm #2}}

\newcommand{\bfb}{{\mathbold b}}

\newcommand{\bfe}{{\mathbold e}}

\newcommand{\bfg}{{\mathbold g}}

\newcommand{\bfl}{{\mathbold l}}

\newcommand{\bfn}{{\mathbold n}}

\newcommand{\bft}{{\mathbold t}}

\newcommand{\bfv}{{\mathbold v}}

\newcommand{\bfx}{{\mathbold x}}

\newcommand{\bfF}{{\mathbold F}}

\newcommand{\bfI}{{\mathbold I}}

\newcommand{\bfL}{{\mathbold L}}

\newcommand{\bfV}{{\mathbold V}}

\newcommand{\bfX}{{\mathbold X}}